\documentclass[a4paper,12pt,draft]{amsart}
\usepackage{amsmath}
\usepackage{amsfonts}
\usepackage{amssymb}
\usepackage[OT4]{fontenc}
\usepackage{fullpage}
\newtheorem{theorem}{Theorem}
\newtheorem{lemma}[theorem]{Lemma}
\newtheorem{corollary}[theorem]{Corollary}
\newtheorem{observation}[theorem]{Observation}

\newtheorem{definition}[theorem]{Definition}

\newtheorem{Question}[theorem]{Question}
\numberwithin{theorem}{section}
\numberwithin{equation}{section}
\DeclareMathOperator\dc{{\it dd^c}}
\DeclareMathOperator\co{{\it\mathbb C^n}}

\DeclareMathOperator\we{\wedge}
\DeclareMathOperator\om{{\it\omega}}

\DeclareMathOperator \caw {{\it cap_{\om}}}
\DeclareMathOperator \wph {{\it\om_{\phi}}}
\DeclareMathOperator \wphj {{\it\om_{\phi_j}}}
\DeclareMathOperator \wps {{\it\om_{\psi}}}
\DeclareMathOperator \wphn {{\it\om_{\phi}^n}}
\DeclareMathOperator \wphnj {{\it\om_{\phi_j}^n}}
\DeclareMathOperator \wpsn {{\it\om_{\psi}^n}}
\DeclareMathOperator \wphnmo {{\it\om_{\phi}^{n-1}}}
\DeclareMathOperator \wpsnmo {{\it\om_{\psi}^{n-1}}}
\DeclareMathOperator \wro {{\it\om_{\rho}}}
\DeclareMathOperator \ex {{\it\mathcal E(X,\om)}}
\DeclareMathOperator \dx {{\it\mathcal D(X,\om)}}
\DeclareMathOperator \dxa {{\it\mathcal D^a(X,\om)}}
\DeclareMathOperator \lbr {{\it\lbrace}}
\DeclareMathOperator \rbr {{\it\rbrace}}

\title{Uniqueness and stability in $\ex$}
\author{S\l awomir Dinew}
\begin{document}

\maketitle

\begin{abstract}  We prove  uniqueness for the Dirichlet problem for the complex Monge-Amp\`ere equation on compact K\"ahler manifolds in the case of  probabilility measures vanishing on pluripolar sets. The proof uses the mass concentration technique due to Ko\l odziej coupled with inequalities for mixed Monge-Amp\`ere measures and the comparison principle. As a by-product we generalize Xing's stability theorem.
\end{abstract}
\section{Introduction}
Pluripotential theory on compact K\"ahler manifolds turned out to be a very effective tool in complex geometry and dynamics. Despite the fact that it deals with (a priori) non smooth functions, the techinques were used with success in purely geometrical problems. For example, the $L^{\infty}$\ estimate from \cite{Koj98}\ gives enough flexibility for studying various limiting problems in geometry which were unaccessible with the standard (restrictive) PDE techniques.

However the ideas considered soon opened new problems of independent interest. In a series of articles Guedj, Zeriahi and collaborators (\cite{GZ1,GZ2,EGZ,BGZ,DGZ}) and Ko\l odziej (\cite{Koj98,Koj03}) laid down the foundations of the theory. In particular the complex Monge-Amp\`ere operator was defined on (non-smooth) $\om$-psh functions (see the next section for  definitions of all the notions appearing in the introduction), and its maximal domain of definition was explorated. By analogy to the works of B\l ocki (\cite{Bl2,Bl3}) and Cegrell (\cite{Ce1,Ce2}) many results form the ''flat'' theory (i.e. the one in domains in $\co$) were adjusted to the K\"ahler manifold case. However from the very beginning some problems revealed to be unexpectedly difficult in the  new setting. One of them is the problem of uniqueness. 

This problem  consists of the following:
Let $\phi,\ \psi$\ be $\om$-psh functions and the Monge-Amp\`ere operator is well defined for them. Assume that $\wphn=\mu=\wpsn$. The question is under what assumptions on the positive measure $\mu$ and/or on the functions $\phi,\ \psi$\ one can conclude that $\phi-\psi$\ is constant?

Note that there are known examples of failure of the uniqueness in general (see \cite{Bl1} for one such example), hence some assumptions are necessary.

The first result in this direction was done by E. Calabi (\cite{Ca}). He proved that if $\phi, \psi$\ are smooth and $\wph,\ \wps$\ are K\"ahler forms (i.e. strictly positive) then uniqueness does hold. These are natural assumptions from geometer's perspective and the proof is quite easy in this case. However both smoothness and strict positivity are crucial in this approach, hence it gives no insight what to do in general.

The next step was done by Bedford and Taylor \cite{BT89} who  proved uniqueness for bounded $\phi, \psi$\ provided the underlying mainfold is $\mathbb P^n$. Their main idea was to control the $L^2$\ norm of the gradient of the difference of $\phi$\ and $\psi$.

Using different technique Ko\l odziej \cite{Koj03} proved uniqueness for bounded functions on arbitrary compact K\"ahler manifold modulo additional mild assumptions on the measure $\mu$.

The ''bounded'' case was finally done by B\l ocki \cite{Bl1}. The proof has some common points with the one in \cite{BT89}, but is much easier and transparent. Furthermore the proof gives some stability results showing that when one perturbs the measure on the right hand side slightly the normalized solution is in a way close to the original one.

By developing theory of Cegrell classes in the K\"ahler manifold setting \cite{GZ1,GZ2} (see also \cite{Di1}) the domain of definition of the operator was enlarged with many unbounded functions. Guedj and Zeriahi \cite{GZ2} observed that B\l ocki's argument, with suitable modifications, can be carried over to prove uniqueness in the class $\mathcal E^1(X,\om)$. Recently Demailly and Pali \cite{DP} proved uniqueness in the same class for more general forms which are only semi-positive. The most general result so far was proven very recently by B\l ocki (see \cite{Bl4}) who proved that uniqueness does hold in the class $\mathcal E^{1-\frac1{2^{n-1}}}(X,\om),\ n=dim X$.

Simultaneously the picture in the flat theory was made much clearer by Cegrell who proved in \cite{Ce2} that one can prove uniqueness provided the measure $\mu$\ does not charge pluripolar sets. The proof however relied heavily on tools that are not available in the K\"ahler setting. Nevertheless it is natural to expect that uniqueness in this class (called $\ex$) should also hold. In fact in \cite{DGZ} this point is an important obstruction for further understanding of the domain of definition of the Monge-Amp\`ere operator.

The class $\ex$\ deserves special interst, due to the following result proven in \cite{GZ2}:
\begin{theorem}
Let $\mu$\ be a probability Borel measure that vanishes on pluripolar sets. Then there exists (at least one) $\rho\in\ex$\ such that 
$$\om_\rho^n=\mu,\ \sup_X\ \rho=0.$$
\end{theorem}
 It is therefore important to study what happens in $\ex\setminus\mathcal E^1$, so one can understand better the action of the complex Monge-Amp\`ere operator. Despite numerous attempts the problem remained open until now. Let us state our main result:
\begin{theorem}
Let $\phi,\ \psi\in\ex,$\ be such that $\wphn=\wpsn$. Then $\phi-\psi$\ is constant.
\end{theorem}\label{mainres}
In \cite{CK} a stability theorem for $\mathcal D^a(\Omega)$\ was proven, where $\Omega$\ is a hyperconvex domain in $\co$. A question appears if such stability does hold in the K\"ahler setting either.

 In the case of bounded functions such stability result was proven in \cite{Koj05}, (there are some immaterial flaws in the proof there). 

 Recently Y. Xing \cite{Xi} proved the corresponding stability result in $\mathcal E^1(X,\om)$. Note that all his proof goes through in $\ex$ too, provided we have the uniqueness result. 
 Hence  it is now straightforward to generalize this in $\ex$:
\begin{theorem}\label{stability}
Suppose $\mu$\ is a measure that vanishes on all pluripolar sets.  Suppose we have a sequence  $\phi_j\in\ex,\ sup_X\phi_j=0\ \wphnj\leq\mu,\ \forall j\geq 1$. Let also $\phi\in\ex,\ sup_X\phi=0$\ solve $\wphn\leq\mu$. Then the following are equivalent:
\begin{enumerate}
\item $\wphnj\rightarrow \wphn$\ as measures,
\item $\phi_j\rightarrow\phi$\ in $L^1$.
\end{enumerate}
Furthermore in this case the functions $\phi_j$\ converge to $\phi$\ in capacity. 
\end{theorem}
 
  Note that convergence in capacity (for $\om$-psh functions) is stronger than convergence in $L^p$\ norm $1\leq p<\infty$, but is weaker than uniform convergence. Since we deal with a priori unbounded and discontinuous functions this is in a way the best convergence one can expect.

{\bf Acknowledgement} The Author wishes to express his sincere gratitude to professors Zbigniew B\l ocki and S\l awomir Ko\l odziej for numerous discussions on the topic.
\section{preliminaries}
Throughout the note we shall work on a fixed compact $n$-dimensional K\"ahler manifold $X$\ equipped with a fundamental K\"ahler form $\omega$\ (that is $d$-closed strictly positive globally defined form) given in local coordinates by 
$$\omega={\frac i 2}  \sum_{k,j=1}^n g_{k\overline{j}}dz^k \wedge d\overline z^j$$
 We assume that the metric is normalized so that 
$$\int_X\omega^n=1.$$ 
Recall that 
\begin{equation*}
 PSH(X,\omega):=\lbrace \phi \in L^1(X,\omega):dd^c\phi \geq -\omega,\ \phi \in\mathcal C^{\uparrow}(X) \rbrace
\end{equation*} 
where $\mathcal C^{\uparrow}(X)$\ denotes the space of upper semicontinuous functions and, as usual, $d$\ is the standard operator of exterior differentiation while $d^c:=i/2\pi(\overline{\partial}-\partial)$. We call the functions that belong to $PSH(X,\omega)$ $\omega$-plurisubharmonic ($\om$-psh for short). We shall often use the handy notation $\wph:=\om+dd^c\phi$.

  These  functions are locally standard plurisubharmonic functions minus a (smooth) potential for the form $\omega$. This allows to use classical local results from pluripotential theory. In particular the Monge-Amp\`ere operator
$$\wphn:=\om_{\phi}\we\cdots\we\om_{\phi}$$
 is well defined  for bounded $\om$-psh functions. This approach was used in \cite{Koj98,Koj03,Di1}.
 
 We recall below the definition of the class $\ex$.

For every $u\in PSH(X,\omega)$\ $(\omega+\dc max( u,-j))^n$\ is a well defined probability measure. By \cite{GZ2} the sequence of measures $\chi_{\lbrace u>-j\rbrace}(\omega+\dc max( u,-j))^n$\ is always increasing and one defines
$$\mathcal E(X,\omega):=\lbrace u\in PSH(X,\omega)\ |\ \lim_{j\rightarrow\infty}\int_{X}\chi_{\lbrace u>-j\rbrace}(\omega+\dc max( u,-j))^n=1\rbrace.$$ 
These functions are a priori unbounded, but the integral assumption ensures that the Monge-Amp\`ere measure has no mass on $\lbrace u=-\infty\rbrace$. Then one defines 
$$(\omega+\dc u)^n:=\lim_{j\rightarrow\infty}\chi_{\lbrace u>-j\rbrace}(\omega+\dc max( u,-j))^n.$$
In particular Monge-Amp\`ere measures of functions from $\mathcal E(X,\omega)$\ do not charge pluripolar sets. We refer to \cite{GZ2} for a discussion of that notion.

The class $\mathcal E^1(X,\om)$, or more generally $\mathcal E^p(X,\om), p>0$\ is defined by
$$\mathcal E^p(X,\om):=\lbr\phi\in\ex\ |\ \int_X|\phi|^p\om_{\phi}^n\leq\infty\rbr.$$
Since $\om$-psh functions are upper semicontiunuous, they are bounded from above, hence one usually considers only nonnegative $\om$- psh functions from $\mathcal E^p(X,\om)$, which often comes in handy in technical details. Note that originally the classes $\mathcal E^p$\ were defined (similarly to the Cegrell classes in the flat theory) with the use of a  sequence of bounded functions $\phi_j, \phi_j\searrow\phi$, such that $\sup_j\int_X |\phi_j|^p\om_{\phi_j}^n\leq\infty$. The results form \cite{GZ2,Di1} have shown that actually one can take just the sequence $\phi_j:=\max (\phi,-j)$, hence both definitons are coherent.

 One can  also define ''local''  classes in an attempt similar to the one from \cite{Bl2,Bl3}.
 We define the class $\dx$\ by
$$\dx:=\lbr\phi\in PSH(X,\om)\ |\ \forall z\in X\ \exists U_z- {\rm open},\ z\in U_z,\ \rho+\phi\in \mathcal D(U_z) \rbr,$$
where $\rho$\ is a local potential in $U_z$\ for $\om$\ and $\mathcal D(U_z)$\ is the maximal domain of definition of the Monge-Amp\`ere operator in $U_z$, (see \cite{Bl2,Bl3}). 

 Note however that the ''local'' and global definition yield different classes, as shown in \cite{GZ2}. This is in sharp contrast with the ''flat'' theory.

Define also $\dxa$\ by
$$\dxa:=\lbr\phi\in\dx\ |\ \wphn(A)=0,\ \forall A\subset X,\ A-{\rm pluripolar}\rbr.$$
It is known that $\dxa\subset\ex$, while $\dx\nsubseteq\ex\nsubseteq\dx$.

Please note that the terminology in the Cegrell classes, partially due to the mentioned differences in the ''local'' and ''global'' settings varies in the literature. In particular the class $\dx$\ is denoted by $\ex$\ in \cite{HKH} or \cite{Di2}. The class $\ex$\ in turn differs in some aspects from the class $\mathcal E(\Omega)$\ in the ''flat'' setting (for example a function in $\mathcal E(\Omega)$\ may have a Monge-Amp\`ere measure that charges points).

The first result that we shall need- an inequality for mixed Monge-Amp\`ere measures was shown in a special case in \cite{Koj03} and in full generality in \cite{Di2}.
\begin{theorem}\label{maineq}
Let $u,\ v\in\mathcal E(X,\omega)$\ be $\omega$-psh functions,  $\mu$\ be a positive measure that does not charge pluripolar sets and $f,\ g\in L^1(d\mu)$. If
$$(\omega+\dc u)^n\geq fd\mu,\ (\omega+\dc v)^n\geq gd\mu$$
as measures, then 
$$(\omega+\dc u)^k\we(\omega +\dc v)^{n-k}\geq f^{\frac kn}g^{\frac{n-k}n}d\mu,\ \forall k\in\lbrace 1,\cdots,n-1\rbrace$$
\end{theorem}
\begin{corollary}\label{corol}
If $\phi,\ \psi\in\ex$\ and $\wphn=\wpsn$, then also for every $ t\in(0,1)$ we have
$$\om_{t\phi+(1-t)\psi}^n=\wphn=\wpsn.$$
\end{corollary}
The next result we shall need is somewhat nonstandard, so although it is similar to the usual comparison principle, we sketch a proof.
\begin{theorem}[''partial'' comparison principle]
Suppose $T$\ is a $(k,k)$\ positive closed current on $X$\ of the form $\om_{\phi_1}\we\cdots\we\om_{\phi_k},\ \phi_j\in\ex\ \forall j\in\lbr1,\cdots,k\rbr$, where $0\leq k\leq n-1$. Let furthermore $u,\ v\in\ex$. Then
$$\int_{\lbr u<v\rbr}\om_v^{n-k}\we T\leq \int_{\lbr u<v\rbr}\om_u^{n-k}\we T.$$
\end{theorem} 
\begin{proof}
Note that in the case $k=0$\ this is the standard comparison principle in $\ex$, which was shown in Theorem 1.5 in\cite{GZ2},  (historically the bounded case was first shown in \cite{Koj03}). Note also that it is enough to get the statement for $n-k=1$, since all the other cases can be done by iteration of the $n-k=1$\ argument. So, we assume $n-k=1$.

In the class $\dxa$\ this was proven in \cite{HKH}. Our case is indeed not much different:

 Basically one can repeat the argument from the Theorem 1.5 in \cite{GZ2}, provided one knows that for bounded $u,\ v$
$$\chi_{\lbr u<v\rbr}\om_v\we T=\chi_{\lbr u<v\rbr}\om_{\max(u,\ v)}\we T,$$ 
corresponding to equality (1) in the proof of Guedj and Zeriahi.
Now one can proceed in the following way: define the canonical approximants for $\phi_1,\cdots,\phi_{n-1}$\ by
$$\phi_s^{(j)}:=max(\phi_s,-j).$$ 
Let $T^{(j)}:=\om_{\phi_1^{(j)}},\we\cdots\we\om_{\phi_{n-1}^{(j)}}$. Now all the functions appearing in the wedge products belong to $\dxa$, hence \cite{HKH} applies. So, for $u,\ v$\ bounded we have
$$\int_{\lbr u<v\rbr}\om_v\we T^{(j)}\leq \int_{\lbr u<v\rbr}\om_u\we T^{(j)}.$$
By propreties of $\ex$\ proven in \cite{GZ2} for bounded $u,\ v$\ passing to the limit with $j\rightarrow\infty$\ one gets
$$\int_{\lbr u<v\rbr}\om_v\we T\leq \int_{\lbr u<v\rbr}\om_u\we T.$$
Now the step from bounded $u,\ v$\ to general $u,\ v$\ can be done exactly as in Theorem 1.5 in \cite{GZ2}.

\end{proof}

{\bf Remark.} It is likely that the result above holds for every positive closed current $T$\ of bidegree $(k,k)$. Also it is likely that the partial comparison principle should hold for more general functions than the $\om$-psh ones and this new class should depend on $T$. We shall not go into details since the stated version is satisfactory for our needs.\\

For our discussion of stability we shall need the notion of Bedford-Taylor capacity and convergence with respect to it, first used in the context of K\"ahler manifolds in \cite{Koj03}:
\begin{definition} Given a Borel set $A\subset X$, we define the Bedford-Taylor capacity as
$$\caw(A):=sup\lbr\int_A\om_u^n\ |\ u\in PSH(X,\om),\ 0<u<1\rbr.$$
\end{definition}
\begin{definition}
We say that a sequence $\phi_j\in PSH(X,\om)$\ converges in capacity to $\phi\in PSH(X,\om)$\ if
$$\forall t>0\ \caw(|\phi_j-\phi|>t)\rightarrow 0,\ j\rightarrow\infty.$$
\end{definition}
We refer to \cite{Koj98}, \cite{GZ1}\ for the basic properties of this capacity and the notion of convergence with respect to it. Recently Hiep in \cite{Hi} obtained the following characterization of convergence in capacity for uniformly bounded functions:
\begin{theorem}\label{Hiep} Let $\phi_j,\ \phi$\ be uniformly bounded $\om$-psh functions the following are equivalent:
\begin{enumerate}
\item $\phi_j$\ converges to $\phi$\ in capacity,
\item $\limsup_{j\rightarrow\infty}\phi_j\leq \phi$\ and $\int_X(\phi_j-\phi)\wphnj\rightarrow 0$.
\end{enumerate}
\end{theorem}

\section{proofs}
\begin{proof}[Proof of Theorem \ref{mainres}]

Consider the level sets
$$A_t:=\lbr \phi-\psi=t\rbr,\ t\in\mathbb R\cup\lbr+\infty\rbr\cup\lbr-\infty\rbr.$$

These are all Borel sets which are closed in the plurifine topology. The main ingredient of the proof is to show that the whole mass of $\mu:=\wphn=\wpsn$\ is concentrated on exactly {\bf one} of the sets $A_t$.
To achieve this we need some additional results

\begin{lemma} The measure charges at most countably many of the sets $A_t$ and does not charge neither $A_{+\infty}$\ nor $A_{-\infty}$.
\end{lemma}

\begin{proof}[Proof.] The first part is proved in \cite{GZ2} Corollary 1.10, and the second follow from the fact that both $A_{+\infty}$\ and $A_{-\infty}$\ are pluripolar.
\end{proof}

Suppose now that the mass of $\mu$\ is not concentrated on one set $A_t$. Then we claim that we can find $t_0\in\mathbb R$\ and a constant $1/2<q<1$\ such that:
\begin{enumerate}
\item $\int_{A_{t_0}}d\mu=0$,
\item $\int_{\lbr\phi<\psi+t_0\rbr}d\mu<q$,
\item $\int_{\lbr\phi>\psi+t_0\rbr}d\mu<q$.
\end{enumerate}
Indeed, One can find $t_1\in\mathbb R$\ such that
$$ 0<\int_{\lbr\phi<\psi+t_1\rbr}d\mu<1,$$
(for othervise the whole mass is concentrated on one level set). Now if $\int_{A_{t_1}}d\mu=0$, we take $t_0:=t_1,\ q=\max\lbr\int_{\lbr\phi<\psi+t_1\rbr}d\mu,\ 1-\int_{\lbr\phi<\psi+t_1\rbr}d\mu\rbr+\epsilon$, with $\epsilon>0$\ so small that still $q<1$. If $A_{t_1}$\ is charged, then for almost every $t<t_1$ $A_t$\ is not charged and by monotone convergence one can take $t_2<t_1$\ close enough to $t_1$\ such that both $A_{t_2}$\ is massless and still $0<\int_{\lbr\phi<\psi+t_1\rbr}d\mu<1$. Take $t_0:=t_2,\ q$\ defined as before and we get the desired properties.

Since adding a constant to $\phi$\ or $\psi$\ is harmless for our discussion we  assume from now on that $t_0=0$.

Consider the new measure 
$$\widehat{\mu}:=\begin{cases} (1/q)\mu, &on \lbr\phi<\psi\rbr\\
  c\mu, &on \lbr\phi\geq\psi\rbr,              \end{cases}$$
where $c$\ is a nonnegative normalization costant so that $\widehat{\mu}$\ is a nonnegative probability measure (note that this is possible, since, by assumption, $\mu$\ charges the set $\lbr\phi\geq\psi\rbr$).

Of course $\widehat{\mu}$\ does not charge pluripolar sets either (and is also a Borel measure since the set $\lbr\phi\geq\psi\rbr$\ is Borel). By \cite{GZ2} we can solve the Monge-Amp\`ere equation
$$\om_\rho^n=\widehat{\mu},\ \rho\in\ex,\ sup_X\rho=0.$$
Note that at this moment we do not know if $\rho$\ is uniquely defined: we just choose one solution.

In such a case we have a set inclusion
$$U_t:=\lbr(1-t)\phi<(1-t)\psi+t\rho\rbr\subset\lbr\phi<\psi\rbr$$
for every $t\in(0,1)$. Hence on $U_t$ we have
$$\wph^{n-1}\we\om_{(1-t)\psi+t\rho}=(1-t)\mu+t\wphnmo
\we\wro\geq(1+((1/q)^{1/n}-1)t)\wphn,$$
where we have made use of Theorem \ref{maineq}\ and Corollary \ref{corol}.

So by the comparison principle 
\begin{align*}
&(1+((1/q)^{1/n}-1)t)\int_{U_t}\wphn\leq\int_{U_t}\wphnmo
\we\om_{(1-t)\psi+t\rho}\leq\\
&\leq\int_{U_t}\wphnmo\we\om_{(1-t)\phi+t0}=(1-t)\int_{U_t}
\wphn+t\int_{U_t}\wphnmo\we\om.
\end{align*}
 Rearranging terms we obtain
\begin{equation}\label{1}
(1/q)^{1/n}\int_{U_t}\wphn\leq\int_{U_t}\wphnmo\we\om.
\end{equation}
Note that exchanging $\wphnmo$\ with $\wpsnmo$\ in the argument above gives
\begin{equation}\label{2}
(1/q)^{1/n}\int_{U_t}\wpsn\leq\int_{U_t}\wpsnmo\we\om.
\end{equation}
(again we make use of Corollary \ref{corol} here).
Now let $t\searrow 0$. The sets $U_t$\ form an increasing sequence, and $U_t\nearrow\lbr\phi<\psi\rbr\setminus\lbr\rho=-\infty\rbr$. But both measures $\wphn$\ and $\wphnmo\we\om$\ do not charge pluripolar sets, hence we obtain
\begin{equation}\label{3}
(1/q)^{1/n}\int_{\lbr\phi<\psi\rbr}\wphn
\leq\int_{\lbr\phi<\psi\rbr}\wphnmo\we\om.
\end{equation}

One can do this reasoning also on the set $\lbr\phi<\psi\rbr$. Namely we find a measure defined like $\widehat{\mu}$, but with respect to the set  $\lbr\phi<\psi\rbr$. Fixing $\wphnmo$\ (or $\wpsnmo$) and arguing the same way we obtain
\begin{equation}\label{4}
(1/q)^{1/n}\int_{\lbr\phi>\psi\rbr}\wphn
\leq\int_{\lbr\phi>\psi\rbr}\wphnmo\we\om.
\end{equation}
But adding these inequalities and the assumption that $A_0$\ is massless one obtains
\begin{align*}
&(1/q)^{1/n}=(1/q)^{1/n}\int_{\lbr\phi>\psi\rbr}\wphn+
(1/q)^{1/n}\int_{\lbr\phi<\psi\rbr}\wphn\leq\\
&\leq \int_{\lbr\phi>\psi\rbr}\wphnmo\we\om+
\int_{\lbr\phi<\psi\rbr}\wphnmo\we\om\leq 1,
\end{align*}
a contradiction.

So we can assume that the whole mass of $\mu$\ is concentrated on $\lbr\phi=\psi\rbr\neq X$.

The next step will be to prove that the same  holds for all the measures $$\om_{\phi}^k\we\om_{\psi}^{n-1-k}\we\om\ k\in\lbr0,\cdots,n-1\rbr.$$
 Note that for these measures one cannot apply straightforwardly the inequality for mixed Monge-Amp\`ere measures, hence we don't have that these measures are (a priori) globally equal. But we shall  need only the fact that they all are massless on $\lbr\phi\neq\psi\rbr$- and this is exactly what the argument will give.

Let $\phi_j:=max\lbr\phi,-j\rbr$. Fix $t\in (0,1)$. Consider the sets
$$V_{t,j}:=\lbr\phi+(t/j)\phi_j+(3/2)t
<\psi\rbr\subset\lbr\phi<\psi\rbr.$$
Again by the comparison principle we obtain
$$\int_{V_{t,j}}\wph^{k}\we\wps^{n-1-k}\we(\wps+(t/j)\om)\leq
\int_{V_{t,j}}\wph^{k}\we\wps^{n-1-k}\we(\wph+(t/j)\om_{\phi_j}).$$
Now, (recall $\wph^{k}\we\wps^{n-k}=\wphn,\ \forall k\in\lbr0,\cdots n-1\rbr$) the equation above reads
$$\int_{V_{t,j}}\wph^{k}\we\wps^{n-1-k}\we\om\leq
\int_{V_{t,j}}\wph^{k}\we\wps^{n-1-k}\we\om_{\phi_j}.$$
Note that $V_{t,j}$\ is a decreasing sequence of sets in terms of $j$. Letting $j\rightarrow\infty$\ and using vanishing on pluripolar sets we obtain
$$\int_{\lbr\phi+(3/2)t<\psi\rbr}\wph^{k}\we\wps^{n-1-k}\we\om\leq
\int_{\lbr\phi+(3/2)t<\psi\rbr}\wphn=0.$$
Finally letting $t\searrow 0$\ we obtain
$$\int_{\lbr\phi<\psi\rbr}\wph^{k}\we\wps^{n-1-k}\we\om\leq
\int_{\lbr\phi<\psi\rbr}\wphn=0.$$
Again exchanging $\phi$\ with $\psi$\ and $\lbr\phi<\psi\rbr$ with $\lbr\phi>\psi\rbr$\ one obtains that the measures  $\wph^{k}\we\wps^{n-1-k}\we\om$\ are massless on $\lbr\phi\neq\psi\rbr$.

 But then exchanging $\wph^{k}\we\wps^{n-1-k}$\ in the argument above with  $\om_{\phi}^k\we\om_{\psi}^{n-2-k}\we\om$\ we obtain
that $\om_{\phi}^k\we\om_{\psi}^{n-2-k}\we\om^2$\ is again massless on $\lbr\phi\neq\psi\rbr$\ (this time we do not have equality for the mixed Monge-Amp\`ere measures in the whole $X$\ but instead we use vanishing on $\lbr\phi\neq\psi\rbr$). 

 Iterating the argument yields that all the measures
$$\wph^k\we\wps^s\we\om^{n-k-s}$$
are massless on $\lbr\phi\neq\psi\rbr$. Hence this holds also for $\om^n$\ i.e.
$$\int_{\lbr\phi\neq\psi\rbr}\om^n=0.$$
But $\lbr\phi\neq\psi\rbr$\ is a pluri fine open set, hence (if it is non empty) it has positive Lebesgue measure. This is a contradiciton.
\end{proof}
\begin{proof}[Proof of Theorem \ref{stability}]
Recall once again that the equivalence $(1)\Leftrightarrow(2)$\ is essentially due to Xing. In particular Corollary 1 in \cite{Xi}\ yields $(2)\Rightarrow(1)$. To show that $(1)\Rightarrow(2)$\ we proceed exactly as in \cite{Xi} Theorem 8: We can extract a subsequence from $\phi_j$\ tending in $L^1$\ to some function $\psi$. It follows that $
\wphn=\wpsn$, so by uniqueness (now in $\ex$) $\phi=\psi$. Thus, since any subsequence has subsequence convergent in $L^1$\ to $\phi$\ we conclude that $\phi_j$\ converges to $u$\ in $L^1$.

It remains to prove convergence in capacity. Regrettably we do  not have an easy proof in the spirit of \cite{CK,Koj05} for this.

Before we start the proof we state two technical lemmata. Their proof will be postponed until the end.
\begin{lemma}\label{phij} We have an uniform estimate for the sublevel sets
$$\forall N>0,\ \int_{\lbr\phi_j<-N\rbr}\wphnj\leq \epsilon(N),$$
where $\epsilon(N)\searrow 0$\ as $N\rightarrow\infty$.
\end{lemma}
\begin{lemma}\label{small} With the notation as above we have that for $k\in\lbr0,\cdots,n\rbr $
\begin{align*}
&a_k)\forall \delta>0\ {\rm we\ have}\  \int_{\lbr\phi>\phi_j+\delta\rbr}\om^k\we\om_{\phi_j}^{n-k}\rightarrow 0,\\
&b_k)\om^k\we\om_{\phi_j}^{n-k}\rightarrow\om^k\we\om_{\phi}^{n-k}\ {\rm as\ measures}.
\end{align*}
\end{lemma}

Indeed, suppose this is already shown. We shall apply the same idea as in \cite{CK} (see also \cite{Koj05}) suitably adapted to K\"ahler setting. Consider the set $\lbr\phi>\phi_j+\delta\rbr$. Fix a function $h\in PSH(X,\om),\ -1<h<0$. Then we have  set inclusions
\begin{equation}
\lbr\phi>\phi_j+\delta\rbr\subset 
\lbr\phi+(\delta/2)h>\phi_j+\delta/2\rbr
\subset\lbr\phi>\phi_j+\delta/2\rbr.
\end{equation}
By the comparison principle
\begin{align*}
&(\delta/2)^n\int_{\lbr\phi>\phi_j+\delta\rbr}\om_h^n\leq
\int_{\lbr\phi+(\delta/2)h>\phi_j+\delta/2\rbr}((1+\delta/2)
\om_{\phi+(\delta/2)h})^n\leq\\
&\leq\int_{\lbr\phi>\phi_j+\delta/2\rbr}
((1+\delta/2)\om_{\phi_j})^n\leq\sum_{k=0}^n\binom nk(\delta/2)^k\int_{\lbr\phi>\phi_j+\delta/2\rbr}
\om^k\we\om_{\phi_j}^{n-k}.
\end{align*}
But, according to Lemma \ref{small}, all the integrals in the right hand side tend to zero as $j$\ increases to $\infty$. So
\begin{equation}
\int_{\lbr\phi>\phi_j+\delta\rbr}\om_h^n\rightarrow0,\ \forall \delta>0.
\end{equation}
Now since $h$\ was arbitrary, it follows that
$$\caw(\lbr\phi>\phi_j+\delta\rbr)\rightarrow 0.$$
It is a consequence of Hartogs' lemma that
$$\caw(\lbr\phi_j>\phi+\delta\rbr)\rightarrow 0.$$
Coupling these we obtain that $\phi_j$\ converges to $\phi$\ in capacity.

So, below we prove Lemmata \ref{phij}\ and \ref{small}.

\begin{proof}[Proof of Lemma \ref{phij}]
Note that contrary to the ''flat setting'' we do not have that 
$$\forall u\in\ex\ \limsup_{s\rightarrow\infty}s^n\caw(\lbr u\leq-s\rbr)=0,$$
(see \cite{GZ2} for capacity estimates of sublevel sets), so here we make essential use of the uniform domination of the measures. Indeed, it is enough to prove that
$$\int_{\lbr\phi_j<-N\rbr}\mu\rightarrow 0.$$
Since $\mu$\ does not charge pluripolar sets, we have that $\mu$\ is equicontinuous with $\caw$\ in the Xing's terminology, see \cite{Xi}. This means that $$\forall\epsilon>0\ \exists\delta>0\ \forall  E\subset X\ \caw(E)\leq \epsilon\Rightarrow \mu(E)\leq\delta.$$ 
All we have to do now is to observe that by Proposition 2.6 in \cite{GZ1} we have
$$\forall u\in PSH(X,\om),\ \caw(\lbr u<-s\rbr)\leq C/s,$$
for some constant independent of $u$, so sets $\lbr\phi_j<-N\rbr$\ have uniformly  small capacity (in terms of $j$), as $N$\ increases.
\end{proof} 
\begin{proof}[Proof of Lemma \ref{small}]
The proof will be inductive. We will prove $(a_0)$\ and we will show the implication $(a_k)\Rightarrow(b_k)$. Finally we prove $(a_k)\cup(b_k)\Rightarrow(a_{k+1})$.

 1.{\it Proof of $(a_0)$}.
Note that by assumption it is enough to prove that
\begin{equation}\label{meassmall}\forall \delta>0\ {\rm we\ have}\  \int_{\lbr\phi>\phi_j+\delta\rbr}\mu\rightarrow 0,
\end{equation}

Note however that (\ref{meassmall}) is a standard consequence of vanishing on pluripolar sets. Indeed, by Lemma 1.4 in \cite{CK} we have
$$\lim_{j\rightarrow\infty}\int_{X}(\max(\phi_j,-M)-\max(\phi,-M))d\mu\rightarrow 0,\ \forall M>0,$$
while, by equicontinuity of $\mu$\ with $\caw$, for $M$\ big enough both $\mu(\lbr\phi_j<-M\rbr)$\ and $\mu(\lbr\phi<-M\rbr)$\ can be made arbitrarily small. This concludes the proof of $(a_0)$.

2. {\it Proof of $(a_k)\Rightarrow (b_k)$.}
 Define as in \cite{Koj98}
$$\rho_s:=max(\phi_{j(s)},\phi-1/s),$$
where $j(s)$\ is choosen so big that
 \begin{equation}\label{last}
\int_{\lbr\phi>\phi_{j(s)}+1/s\rbr}\om^k
\we\om_{\phi_{j(s)}}^{n-k}\leq 1/s.
\end{equation}
 For any $N>0$ we have that $\max (\rho_s,-N)$\ converges in capacity to $\max(\phi,-N)$ (by Hartogs' lemma). So we have 
$$\lim_{s\rightarrow\infty}\om^k\we\om_{\rho_s}^{n-k}\rightarrow
\chi_{\lbr\phi>-\infty\rbr}\om^k\we\wph^{n-k}.$$
But exactly as in Proposition 3.2.1 in \cite{Koj98} any accumulation point of $\om^k\we\om_{\phi(j)}^{n-k}$\ would have to be equal to $\om^k\we\wph^{n-k}$\ (to be precise, using  (\ref{last}) and the convergence for $\om^k\we\om_{\rho_s}^{n-k}$ one shows that the accumulation points should be $\geq\om^k\we\wph^{n-k}$, but since both are probability measures equality must hold).

Hence we conclude that
$$\om^k\we\om_{\phi_j}^{n-k}\rightarrow\om^k\we\wph^{n-k}$$
in the sense of measures.

3. {\it Proof of $(a_k)\cup(b_k)\Rightarrow (a_{k+1})$.}
Note first that by the partial comparison principle
$$0\leq\int_{\lbr\phi_j+\delta<\phi\rbr}\om^k\we\wph\we\wphj^{n-k-1}
\leq\int_{\lbr\phi_j+\delta<\phi\rbr}\om^k\we\wphj^{n-k}\rightarrow 0.$$

Consider the sets $$W_{j,s}:=\lbr\phi_j+(\delta/s)\max(\phi_j,-s)+(3/2)\delta<\phi\rbr,\ s>0.$$
Note that
$$\lbr\phi_j+(3/2)\delta<\phi\rbr\subset W_{j,s}\subset\lbr\phi_j+(1/2)\delta<\phi\rbr.$$
By the partial comparison principle we have
\begin{align*}&(\delta/s)\int_
{W_{j,s}}\om^{k+1}
\we\wphj^{n-k-1}\leq\\
&\leq(\delta/s)\int_
{W_{j,s}}\om^k\we
\om_{\max(\phi_j,-s)}\we\wphj^{n-k-1}+\int_
{W_{j,s}}
\om^k\we(\wphj-\wph)\we\wphj^{n-k-1}).
\end{align*}
But the second integral can be made arbitrarily small (independently of $s$) provided $j$\ is big enough, by $(a_k)$\ and what we just proved. So we take $j=j(s)$\  big enough, so that the integral is estimated by $\epsilon\delta/s$\ for some fixed small $\epsilon$. This shows that
\begin{align*}&\int_
{W_{j,s}}\om^{k+1}
\we\wphj^{n-k-1}\leq\int_
{W_{j,s}}\om^k\we
\om_{\max(\phi_j,-s)}\we\wphj^{n-k-1}+\epsilon\leq\\
&\leq \int_
{\lbr\phi_j\leq-s\rbr}\om^k\we
\om_{\max(\phi_j,-s)}\we\wphj^{n-k-1}+\int_
{W_{j,s}}\om^k\we\wphj^{n-k}+\epsilon=\\
&=\int_
{\lbr\phi_j\leq-s\rbr}\om^k\we\wphj^{n-k}+\int_
{W_{j,s}}\om^k\we\wphj^{n-k}+\epsilon.
\end{align*}
Now, in turn, taking $s$\ big enough the first integral is small (independently of $j=j(s)$), while the second is small by assuming $(a_k)$ (recall that $\lbr\phi_j+3/2\delta<\phi\rbr\subset W_{j,s}\subset\lbr\phi_j+1/2\delta<\phi\rbr$). This concludes the proof of $(a_{k+1})$.
\end{proof}
\end{proof}
\section{Observations and Corollaries}
\begin{observation} In the proof of uniqueness we haven't made use of the full strenght of the K\"ahler condition on $\om$. Note that if $\om$\ is merely {\bf big} i.e. $\om\geq0,\ \int_X\om^n>0$, the same argument will go through, provided the set $\om^n=0$\ is of zero Lebesgue measure (for a discussion of the case of semi-positive forms we refer to \cite{BGZ,EGZ,DP}). The proof carries through even if we allow merely continuous (instead of smooth) potentials for $\om$, but the final argument cannot be done without any assumption on the measure of the degeneration set. 
\end{observation}
\begin{observation}
The assumptions from the observation above are exactly the same as those from \cite{DP}. So the result is also a generalization of the Demailly-Pali uniqueness in the big case. Note also that the assumptions are satisfied in many cases that are of interest for dynamics and geometry: for example if $\om$\ is a pull-back of a K\"ahler form via proper holomorphic mapping, the degeneration set is analytic, hence of zero Lebesgue measure.
\end{observation}
\begin{observation} The proof of convergence in capacity could be essentially simplified, if we knew the following:

 Let $\phi_j,\ \phi\in\ex$, $\phi_j$\ converge to $\phi$\ as distributions and $\wphnj\rightarrow\wphn$\ weakly. Is that true that for every $ M>0$ we have
$$\om_{\max(\phi_j,-M)}^n\ {\rm\ tends\ weakly\ to}\ \om_{\max(\phi,-M)}^n?$$
Indeed, if true, this result would allow to reduce the whole argument to uniformly bounded functions $\max(\phi_j,-M),\ \max(\phi,-M)$, (for $M$\ sufficiently big) and then to use Hiep's characterization of convergence in capacity, see Theorem \ref{Hiep}. Note that with the assumption of uniform domination of the Monge-Amp\`ere measures this fact is an easy consequence of the convergence in capacity, but we are unable to show this unconditionally. 
\end{observation}
\begin{corollary} As observed by B\l ocki (see \cite{Bl4}), since $\dxa\subsetneq\ex$\ it follows that given a probability measure $\mu$\ on $X$\ the Dirichlet problem 
$$\wphn=\mu,\ \sup_X\phi=0,\ \phi\in\dxa$$
need not have a solution in general. This is in sharp contrast with the Dirichlet problem in the bounded case, since it is proven in \cite{Koj05}, that being a Monge-Amp\`ere measure of a bounded $\om$-psh function is a local property.

\end{corollary}

S\l awomir Dinew\\
Institute of Mathematics\\
Jagiellonian University\\
ul. Reymonta 4\\
30-059 Krak\'ow\\
Poland\\
\tt slawomir.dinew@im.uj.edu.pl
\end{document}